\documentclass[11pt]{amsart}
\usepackage{stmaryrd}
\usepackage{mathrsfs}
\usepackage{amssymb}

\usepackage{titletoc}
\usepackage{tikz-cd}

\input xy
  \xyoption{all}
\usepackage{amscd}
\usepackage{amsmath, amssymb}
\usepackage{amsfonts}
\usepackage[colorlinks,linkcolor=blue,citecolor=blue, pdfstartview=FitH]{hyperref}

\numberwithin{equation}{section}

\theoremstyle{plain}
\newtheorem{thm}{Theorem}[section]
\newtheorem{theorem}[thm]{Theorem}
\newtheorem{lemma}[thm]{Lemma}
\newtheorem{corollary}[thm]{Corollary}
\newtheorem{proposition}[thm]{Proposition}

\theoremstyle{definition}

\newtheorem{remark}[thm]{Remark}

\newtheorem{definition}[thm]{Definition}

\newtheorem{example}[thm]{Example}

\newtheorem{defn-thm}[thm]{Definition-Theorem}

\DeclareMathOperator{\Pic}{Pic}

\newcommand{\C}{{\mathbb C}}

\renewcommand{\P}{{\mathbb P}}
\newcommand{\Q}{{\mathbb Q}}

\newcommand{\Z}{{\mathbb Z}}

\newcommand{\Hom}{{ Hom}}

\DeclareMathOperator{\CH}{CH}
\DeclareMathOperator{\sgn}{sgn}

\DeclareMathOperator{\rk}{rk}

\newcommand{\btheorem}{\begin{theorem}}
\newcommand{\etheorem}{\end{theorem}}
\newcommand{\bproposition}{\begin{proposition}}
\newcommand{\eproposition}{\end{proposition}}
\newcommand{\bdefinition}{\begin{definition}}
\newcommand{\edefinition}{\end{definition}}
\newcommand{\bcorollary}{\begin{corollary}}
\newcommand{\ecorollary}{\end{corollary}}
\newcommand{\bproof}{\begin{proof}}
\newcommand{\eproof}{\end{proof}}
\newcommand{\bremark}{\begin{remark}}
\newcommand{\eremark}{\end{remark}}
\newcommand{\eexample}{\end{example}}
\newcommand{\bexample}{\begin{example}}

\newcommand{\elemma}{\end{lemma}}
\newcommand{\blemma}{\begin{lemma}}

\renewcommand{\phi}{\varphi}

\newcommand{\ee}{\end{eqnarray*}}
\newcommand{\be}{\begin{eqnarray*}}

\newcommand{\beq}{\begin{equation}}
\newcommand{\eeq}{\end{equation}}

\newcommand{\bd}{\begin{enumerate}}
\newcommand{\ed}{\end{enumerate}}

\renewcommand{\Hom}{\text{Hom}}

\renewcommand{\>}{\rightarrow}

\usepackage{fancyhdr}
\begin{document}
\title{Categorical characterization of quadrics}
\makeatletter
\let\uppercasenonmath\@gobble
\let\MakeUppercase\relax
\let\scshape\relax
\makeatother
\author{ Duo Li}

\maketitle
\begin{abstract}
We give a characterization of smooth quadrics in terms of the existence of full exceptional collections of certain type, which generalizes a result of C.Vial for projective spaces.\end{abstract}

\setcounter{tocdepth}{1}
\section{Introduction}
Let $X$ be a complex variety and $D^b(X)$ be its bounded derived category.
 An object $E\in \mathcal D^b(X)$  is called exceptional if we have  $$\Hom(E,E[l])=\left\{
\begin{aligned}
&\C  \quad\text{if} \quad l=0 \\
 &0  \quad \text{if} \quad l\ne 0
\end{aligned}
\right.$$
An exceptional collection is a sequence $\{E_1\cdots,E_n\}$ of exceptional objects which satisfy that $\Hom(E_j,E_i[l])$ vanishes for all $j>i$ and all $l\in \Z.$  An exceptional collection $\{E_1\cdots,E_n\}$  is full if $\mathcal D$ is generated by $\{E_i\}.$

In general,  full exceptional collections which  consist of coherent sheaves are rare. The existence of such collections would impose strong restrictions on the geometry of the variety.  We now list some  varieties with such collections. \begin{itemize}  \item
 projective spaces:  for any $a\in\Z,$  $\{\mathcal O(a),\mathcal O(a+1),\cdots,\mathcal O(a+n)\}$ is a full exceptional collection of $\P^n$ (see \cite{B}).\\
\item smooth quadric $Q^n\subset \P^{n+1}.$ Kapranov \cite{Ka} shows that \begin{itemize}

\item if $n$ is odd, for any $a\in\Z,$ $\{S, \mathcal O(a),\mathcal O(a+1),\cdots,\mathcal O(a+n-1)\}$ is a full exceptional collection of $Q^n,$
\item if $n$ is even, for any $a\in\Z,$ $\{S^-,S^+, \mathcal O(a),\mathcal O(a+1),\cdots,\mathcal O(a+n-1)\}$ is a full exceptional collection of $Q^n$
\end{itemize} where $S,$ $S^-,$  $S^+$ are certain spinor bundles. \end{itemize} A\emph{  minifold} is a smooth projective variety $X$  whose $D^b(X)$  admits a full exceptional collection $\mathcal C$ of minimal possible  length $\dim X+1.$
Galkin, Katzarkov, Mellit and Shinder classify  \emph{minifolds} up to dimension $4$ (see \cite{GK}).
If $\mathcal C$ purely consists of line bundles,  C. Vial proves  that $X$ is necessarily a projective space (see \cite[Theorem 1.2]{V1}). In \cite{KO}, S. Kobayashi  and   T. Ochiai show that if there exists an ample line bundle $H$ on a smooth projective variety  $X$ with $\dim(X) c_1(H)\le c_1(X),$ then $X$ is isomorphic to a projective space or a smooth quadric. The purpose of this article is to generalize C.Vial's theorem and  prove a categorical  analog of S. Kobayashi  and   T. Ochiai's classification. In view of Kapranov's result, it is reasonable to consider full exceptional collections of length $\dim X+1$ or $\dim X+2.$
\btheorem \label{odd} Let $X$ be a smooth projective variety of  dimension $n$  $(n\ge 3)$. Suppose that there is a full exceptional collection $\mathcal C$ of $D^b(X)$ which consists of coherent sheaves and is of  length $n+1.$ If $\mathcal C$ contains a sub-collection \{$\mathcal L_1,$ $\mathcal L_2,\cdots, $ $\mathcal L_n$\} where $\mathcal L_i$ $(1\le i\le n)$ are line bundles, then $X$ is isomorphic to $\P^n$ or $Q^n.$ Moreover, if $n$ is an even number, $X$ is isomorphic to $\P^n.$
\etheorem

\btheorem \label{even}

Let $X$ be a smooth projective variety of  dimension $n$  ($n\ge 3 $). Suppose that there is a full exceptional collection $\mathcal C$ of $D^b(X)$ which consists of coherent sheaves and is of  length $n+2.$ If   $\mathcal C$ contains a sub-collection \{$\mathcal L_1,$ $\mathcal L_2,\cdots, $ $\mathcal L_n$\} where $\mathcal L_i$ $(1\le i\le n)$ are line bundles, then $X$ is an even-dimensional variety and is   isomorphic  to $Q^n.$
\etheorem
To prove Theorem \ref{odd} and Theorem \ref{even}, we use the method from C. Vial's article \cite{V1} and some technical input (see Lemma \ref{sequence}).\\

\textbf{Convention:}  In this article, a variety is an integral scheme of finite type  defined over $\C$.
\section{A technical lemma}
Firstly, we note an easy fact about series of numbers.
\blemma\label{sequence}
Let $\{a_1,\cdots,a_n\}$ be $n$ distinct integers with $n\ge5.$ Assume that the cardinality of the set $A=\{a_j-a_i|1\le i<j\le n\}$ is $n.$ Then there exists an integer $d$ such that $(a_1,\cdots, a_n)$ is one of the following series:
\bd
\item $(a_1,\cdots,\widehat{a_1+kd},\cdots,a_1+nd)$ where $1\le k \le n-1$ is an integer.
\item $(a_1,a_1+2d,\cdots,a_1+(n-1)d,a_1+(n+1)d)$
\item $\sigma (a_1,\cdots,a_1+(n-1)d)$ where $\sigma$ is  a composition of disjoint permutations $(i_1,i_1+1)\cdots(i_l,i_l+1)$ for some $1\le l \le n.$
\ed
\elemma
\bproof
 We firstly assume  $a_1<a_2<\cdots <a_n$ and claim that $(a_1,\cdots, a_n)$ is  of type $(1)$ and $(2).$ It is easy to verify our claim when $n$ is $5.$ For the case $n>5,$ the cardinality of   $B:=\{a_j-a_i|1\le i<j\le n-1\}$ is $n-1$ or $n-2.$ If $\sharp B$ is $n-2,$ then $(a_1,\cdots,a_{n-1})$ is an arithmetic progression which can be  denoted  by $(a_1,a_1+d,\cdots, a_1+(n-2)d)$ for some $d\in \Z.$ It is easy to see  $a_n$ is $a_1+n d.$ If $\sharp B$ is $n-1,$ by induction, we have $(a_1,\cdots,a_{n-1})=(a_1,\cdots,\widehat{a_1+kd},\cdots,a_1+(n-1)d)$ ($1\le k \le n-2$) or $(a_1,\cdots,a_{n-1})=(a_1,a_1+2d,\cdots,a_1+(n-2)d,a_1+n d).$ In the first case, if $k$ is $1,$ then $a_n$ equals to $a_1+nd$ or $a_1+(n+1)d.$ If $k$ is bigger than $2,$ then $a_n$ equals to $a_1+nd.$ In the second case, we can assume  $a_n=a_1+md$ for some integer $m>n.$ If $m$ is $n+1$ or $n+2,$ then $(n-1)d$ belongs to $A.$ If we have $m>n+2,$ then $(m-2)d$ belongs to $A.$ Both situations lead to the fact  $\sharp A\ge n+1,$  so we can exclude the case $(a_1,\cdots,a_{n-1})=(a_1,a_1+2d,\cdots,a_1+(n-2)d,a_1+n d)$ and prove our claim.

For an arbitrary series $(a_1,a_2,\cdots,a_n),$  there is a permutation $\alpha$ satisfying  $a_{\alpha(1)}<a_{\alpha(2)}<\cdots<a_{\alpha(n)}.$ Note that the cardinality of the set $C:= \{a_{\alpha(j)}-a_{\alpha(i)}|1\le i<j\le n\}$ is less than $\sharp \{a_j-a_i|1\le i<j\le n\}.$ It follows that $\sharp C$ is $n-1$ or $n.$
If $\sharp C$ is $n-1,$ then there exists an integer $d$ such that $(a_{\alpha(1)},\cdots,a_{\alpha(n)})$ is arithmetic and $C$ is $\{d,2d,\cdots, (n-1)d\}.$ Note that for any integer $k\ne 1,$ if $kd$ belongs to $A,$  then  $\sgn(k)\cdot d$ belongs to $A.$  It follows that $\forall k>1,$ $kd\in A$ implies $-kd\notin A.$ Now we can assume $2d\in B,$ then $(a_{\alpha(1)},a_{\alpha(1)}+2d,a_{\alpha(1)}+4d,\cdots),$ as well as ($a_{\alpha(1)}+d,a_{\alpha(1)}+3d,a_{\alpha(1)}+5d,\cdots$), is a subseries of $(a_1,\cdots,a_n).$ If $a_{\alpha(1)}+d$ appears on the right of $a_{\alpha(1)}$ in $(a_1,\cdots,a_n),$ then we have $A=\{-d,d,2d,\cdots,(n-1)d\}.$ If $a_{\alpha(1)}+d$ appears on the left of $a_{\alpha(1)}$ in $(a_1,\cdots,a_n),$ then we have $-d\in A$ and $a_{\alpha(1)}+3d$  appears on the right of $a_{\alpha(1)}.$  It follows that $A$ is $\{-d,d,2d,\cdots,(n-1)d\}.$ So if $\sharp C$ is $n-1,$ the series $(a_1,\cdots,a_n)$ is of type $(3).$
Now suppose that $\sharp C$ is $n,$ then there exists an  integer $d$ such that $(a_{\alpha(1)},\cdots,a_{\alpha(n)})$ is of type $(1)$ or $(2).$ For any integer $l>0,$ the numbers $ld$ and $-ld$ can not belong to $A$ at the same time. So we can assume $d\in A$ and $-d\notin A.$ By a similar argument, one can  verify $(a_1,\cdots,a_n)=(a_{\alpha(1)},\cdots,a_{\alpha(n)}).$
\eproof
\bremark\label{4}

For the case $n=4,$ $(a_1,a_2,a_3,a_4)$ is one of the following series:
\bd
\item
$(a,a+d,a+2d,a+4d)$ or $(a,a+2d,a+3d,a+4d)$ where $d$ is an integer.
\item
$\sigma (a,a+d,a+2d,a+3d)$ where $d$ is a positive integer and $\sigma\ne id$ is a permutation.
\item
$\sigma (a,a+d,b,b+d)$ where  $d$ is a positive  integer with $b>a+d$ and $\sigma$ is a permutation.
\ed
\eremark
\section{Proof of Theorem \ref{odd} and Theorem \ref{even}} \noindent\\
Let $X$ be a smooth projective variety. Recall that, for any two objects $E$ and $F$ in $D^b(X),$ the Euler pairing $\chi$ is the integer $\chi(E,F):=\sum_l(-1)^l \dim_{\C}\Hom(E,F[l]).$
\bdefinition
An object $E$ is said to be numerically exceptional if $\chi(E,E)=1.$ A collection $\{E_1,\cdots, E_r\}$ of numerical exceptional objects is a numerical exceptional collection if $\chi(E_j,E_i)=0$ for all $j>i.$ \edefinition
For Fano varieties of Picard number $1,$ we have the following result.
\btheorem\label{M0}
Let Let $X$ be a smooth projective Fano variety of dimension $n$ $(n\ge 3)$ whose $\Pic(X)$ is isomorphic to $\Z.$ Assume that there is a numerical  exceptional collection $\{\mathcal L_1,\cdots \mathcal L_n\}$ where $\mathcal L_i$ $(1\le i \le n)$ are line bundles. Then $X$ is isomorphic to $\P^n$ or a smooth quadric $Q^n.$

\etheorem

 \bproof  We denote $c_1(X)$ by $\lambda H$ for some  $\lambda\in \Z^+.$   The Euler character $\chi(\mathcal O_X(aH))$ is a polynomial $P$ with rational coefficients of degree $n$ in the variable $a.$ We write $\mathcal L_i$ as $\mathcal L_i=O_X(a_iH)$ for some integer $a_i$ and we have $a_i\ne a_j$ for any $i\ne j,$ as $\chi(\mathcal O_X)$ is $1.$  The equalities $\chi(\mathcal L_j, \mathcal L_i)=\chi(\mathcal L_i\otimes \mathcal L_j^{-1})=P(a_i-a_j)=0$  $(i<j)$ imply that $a_i-a_j$ $(i<j)$ are roots of $P(a)=0.$  So the cardinality $\mu$ of the set \{$a_j-a_i| 1\le i<j\le n$\} is $n-1$ or $n.$

If $\mu$ is $n-1,$  the progression $(a_1, a_2,\cdots, a_n )$ is an arithmetic progression  with $a_2-a_1=d\in \Z.$ The polynomial $P$ vanishes at $-ld$ for $1\le l \le n-1.$ Let the remaining root of $P(a)=0$ be $-N.$ By Riemann-Roch, we have $$P(a)=\chi(\mathcal O_X(aH))=\frac{\deg(H^n)}{n!}a^n+\frac{\deg(H^{n-1}\cdot c_1(X))}{2(n-1)!}a^{n-1}+\cdots+\chi(\mathcal O_X).$$
Note  the following  equalities: $$\frac{\deg(H^n)}{n!}(n-1)!d^{n-1}N=1\qquad(1)$$  $$\frac{\deg(H^n)}{n!}(N+\sum_{l=1}^{n-1}ld)=\frac{\deg(H^{n-1}\cdot c_1(X))}{2(n-1)!}\qquad(2).$$
One can deduce    $d=1,$ $\deg H^n=1,$ $\lambda=n+1$ or $d=1,$ $\deg H^n=2,$ $\lambda=n.$ Then by \cite{KO}, the Fano variety $X$ is isomorphic to $\P^n$ or $Q^n.$

If $\mu$ is $n,$  we aim to show that $X$ is isomorphic to $\P^n.$ We firstly assume  $n\ge 5,$ then the progression $(a_1, a_2,\cdots, a_n )$ is classified in Lemma \ref{sequence}. If  $(a_1, a_2,\cdots, a_n )$ is of type $(1)$ as in Lemma \ref{sequence},  then there exists an integer $d$ such that  the polynomial $P(a)$ vanishes exactly at $-ld$ $(1\le l\le n).$ By a similar argument, we have  $\deg H^n=1,$ $d=1$ and $\lambda=n+1.$   If  $(a_1, a_2,\cdots, a_n )$ is of type $(2)$ as in Lemma \ref{sequence}, the roots of  $P(a)=0$ is $\{-d,-2d,\cdots,-(n-1)d,-(n+1)d\}$ for some integer $d.$ Then we have $\frac{\deg(H^n)}{n!}(n-1)!d^{n}(n+1)=1,$ which is absurd as $d$ is an integer. If $(a_1, a_2,\cdots, a_n )$ is of type $(3)$ as in Lemma \ref{sequence},  the roots of  $P(a)=0$ is $\{\pm d,-2d,\cdots,-(n-1)d\}$ for some integer $d.$ Then it is easy to deduce  $\deg(H^n)d^n=-n$ and $\lambda=\frac{(n+1)(n-2)}{n}d\in \Z^+,$ which is impossible.
 Now suppose $n=3.$ The progression of integers $(a_1, a_2, a_3 )$ is not arithmetic. By a simple calculation, we have $\lambda=\frac{4}{3}(a_3-a_1)\ge 4 \in \Z^+.$
Now suppose $n=4.$ If $(a_1,a_2,a_3,a_4)$ is of type $(1)$ as in Remark \ref{4},  by a similar argument, we can deduce that $\lambda$ is $5.$ If $(a_1,a_2,a_3,a_4)$ is of type $(2)$ as in  Remark \ref{4}, then there exists an integer $d$ such that the roots of $P(a)=0$ are $\{d,-d,-2d,-3d\},$
 so we have $\frac{\deg(H^4)}{4!}\prod_{i<j}(a_j-a_i)<0,$ which is impossible. If $(a_1,a_2,a_3,a_4)$ is of type $(3)$ as in  Remark \ref{4}, by the method of exhaustion,  the only possibility satisfying  $\frac{\deg(H^4)}{4!}\prod_{i<j}(a_j-a_i)=1$  and $\lambda>0$ is $(a_1,a_2,a_3,a_4)=(a_1,a_1+1,a_1+3,a_1+4),$  hence $\lambda$ is $5.$
 By all the above argument, if $\mu$ is $n,$ then $\lambda$ is $n+1,$ hence $X$ is isomorphic to $\P^n.$
\eproof
We now prove a lemma about the Picard group of a variety. For similar results, we refer to \cite[Theorem 3.4]{B-P} and \cite[Lemma 2.6]{V1}.
\blemma \label{Pic}
Let $X$ be a smooth projective variety of dimension $n$. Suppose that there is a full exceptional collection $\mathcal C$ of $D^b(X).$  If the length of $\mathcal C$  is $n+1,$  then $\Pic(X)$ is isomorphic to $\Z.$ If the length of $\mathcal C$  is $n+2$ ($n\ge 3$),  then $n$ is an even number and the  $\Pic(X)$ is isomorphic to $\Z.$
\elemma
\bproof For any smooth projective variety $X,$
the existence of a  full exceptional collection of $D^b(X)$ whose length is $m$     implies that the Grothendieck group $K_0(X)$ of $D^b(X)$ is isomorphic to $\Z^m.$ Moreover, we have the identity $$m=\rk K_0(X)=\sum_{i=0}^{n}\rk \CH^i_{\Q}(X)$$ where every $\rk \CH^i_{\Q}(X)$ is  non-zero. If $m$ equals to $n+1,$ then $\Pic(X)$ is of rank $1.$ If $m$ equals to $n+2,$ by \cite{K} or \cite[Proposition 3.10]{V}, the cycle class maps $cl_i:\CH^i_{\Q}(X)\> H^{2i}(X,\Q)$ $(0\le i\le n)$ are all isomorphisms. So for any $0\le i\le n,$ there exists a   perfect pairing $\CH^i_{\Q}(X)\times \CH^{n-i}_{\Q}(X)\>\Q.$ We have $\rk\CH^i_{\Q}(X)=\rk \CH^{n-i}_{\Q}(X)=1$ for $i\ne n-i.$  It follows that $n$ is an even number and $X$ is of Picard number $1.$ By \cite[Lemma 2.2]{G} or \cite[Lemma 2.6]{V1}, the  Chow group $\CH^1(X)$ is torsion free, so $\Pic(X)$ is isomorphic to $\Z.$
\eproof
\noindent \emph{The proof of Theorem \ref{odd}.}
By \cite[Proposition 3.2]{B-P}, the variety $X$ is Fano. Then Theorem \ref{odd} is a direct corollary of Theorem \ref{M0} and Lemma \ref{Pic}.
\qed

\vskip 1\baselineskip
To prove Theorem \ref{even}, we need the concept of the anticanonical pseudoheight which is introduced by A. Kuznetsov.
For an exceptional collection $\mathcal C=\{E_0,\cdots, E_r\},$  the anticanonical pseudoheight is defined as $$ph_{ac}(\mathcal C):=\mathop{\min}\limits_{1\le a_o<\cdots<a_p\le n}(e(E_{a_0},E_{a_1})+\cdots+e(E_{a_{p-1}},E_{a_p})+e(E_{a_p},E_{a_0}\otimes K_X^{-1})-p)$$
where $\forall F,$ $F'\in D^b(X),$  the relative height $e(F,F')$ is defined by $\min\{k|Ext^k(F,F')\ne 0\}.$   A. Kuznetsov shows in \cite[Corollary 6.2]{Ku} that if $ph_{ac}(\mathcal C)>-\dim X,$ then $\mathcal C$ is not full.\\

\noindent \emph{The proof of Theorem \ref{even}.}
By Lemma \ref{Pic}, we have $\Pic(X)\simeq \Z H$ for some ample line bundle $H$ and $\rk K_0(X)$ is $n+2,$ which implies that $X$ is even-dimensional. We keep using the notations in the proof of Theorem \ref{M0} and aim to show that $X$ is Fano.
If $n>5,$ we can list all the possibilities in the following
diagram:

$$\begin{tabular}{|c|c|c|c|c|}
  \hline
  &$(a_1,\cdots,a_n)$ & roots of $P(a)$ & $\lambda$ & $X $\\
  \hline
  (1)& $(a_1,a_1+1,\cdots,a_1+(n-1))$ & $\{-1,-2,\cdots,-(n-1)\}$ & $n$ & $Q^n$\\
  \hline
(2)& $(a_1,a_1+1,\cdots,a_1+(n-1))$ & $\{-1,-2,\cdots,-n\}$ & $n+1$ & $\P^n$\\
  \hline
 (3)&$(a_1,a_1-1,\cdots,a_1-(n-1))$ &  $\{1,2,\cdots,(n-1)\}$ & $-n$ & general type\\
  \hline
  (4)&$(a_1,a_1-1,\cdots,a_1-(n-1))$ & $\{1,2,\cdots,n\}$ & $-n-1$ & general type \\
  \hline

 (5)&$(a_1,\cdots,\widehat{a_1-k},\cdots,a_1-n)$ & $\{1,2,\cdots,n\}$ & $-n-1$ & general type\\
  \hline
\end{tabular}$$

 By the Kodaira vanishing theorem,  for line bundles $\mathcal L,$ $\mathcal L'$ with $\mathcal L\otimes \mathcal L'^{-1}$ ample, the relative height $e(\mathcal L, \mathcal L')$ is bigger than $n.$ Note that for all coherent sheaves $F,$ $F'$, $e(F,F')$ is positive.  It easily follows that in cases $(3)$ $(4)$ $(5)$ of our diagram, $ph_{ac}(\mathcal C)>-n,$ which is a contradiction.

If  $n$ is $4,$ for the cases $\mu=3$ or $\mu=4,$ $(a_1,a_2,a_3,a_4)$ of type $(1)$ and $(2)$ as in Remark \ref{4}, we can easily deduce that $X$ is Fano, otherwise there is $ph_{ac}(\mathcal C)>-4.$ If $(a_1,a_2,a_3,a_4)$ is of type $(3).$ Note that if $\chi(\mathcal L_j, \mathcal L_i)$ vanishes, by the Serre duality, $\chi(\mathcal L_i, \mathcal L_j\otimes K_X)$ also vanishes. It follows that for any $i<j,$ $a_j-a_i-\lambda$ is also a root of $P(a)=0.$ Then by the method of exhaustion, $(a_1,a_2,a_3,a_4)$ is $(a,a+1,a+3,a+4),$ which also implies that $X$ is also Fano. Note that $K_0(\P^n)$ is isomorphic to $\Z^{n+1},$ by Theorem \ref{M0}, $X$ is isomorphic to $Q^n.$

\vskip 1\baselineskip

\textbf{Acknowledgments.}  The author is  very grateful to Professor Baohua
Fu, Yihua Chen,  Zili Zhang for their helpful suggestions and discussions. The author wishes to thank Professor A. Kuznetsov, C. Vial and Shizhuo Zhang for many valuable suggestions for the earlier version of this article.
\end{document}